\documentclass[12pt,reqno]{amsart}

\DeclareMathOperator{\ann}{ann}%
\DeclareMathOperator{\can}{can}%
\DeclareMathOperator{\chr}{char}%
\DeclareMathOperator{\Gal}{Gal}%
\DeclareMathOperator{\res}{res}%
\DeclareMathOperator{\supp}{supp}%
\DeclareMathOperator{\trdeg}{tr.deg}%

\usepackage{amssymb,latexsym,amsmath,amscd}
\usepackage[all]{xy}

\begin{document}

\newtheorem{theorem}{Theorem}
\newtheorem{corollary}{Corollary}
\newtheorem{lemma}{Lemma}
\newtheorem{proposition}{Proposition}

\theoremstyle{definition}
\newtheorem*{remark*}{Remark}

\newcommand{\C}{\mathbb{C}}
\newcommand{\F}{\mathbb{F}}
\newcommand{\Fp}{\mathbb{F}_p}
\newcommand{\Jc}{\mathcal{J}}
\newcommand{\Kf}{\mathfrak{K}}
\newcommand{\N}{\mathbb{N}}
\newcommand{\Q}{\mathbb{Q}}
\newcommand{\R}{\mathbb{R}}
\newcommand{\Z}{\mathbb{Z}}

\newcommand{\mo}{\mathopen\langle}
\newcommand{\mc}{\mathclose\rangle}

\parskip=10pt plus 2pt minus 2pt

\title[Galois Module Construction and Classification]
{Construction and Classification of Some Galois Modules}

\author[J\'{a}n Min\'{a}\v{c}]{J\'an Min\'a\v{c}$^{*\dagger}$}
\address{Department of Mathematics, Middlesex College,
\ University of Western Ontario, London, Ontario \ N6A 5B7 \
CANADA}
\thanks{$^*$Research supported in part by the Natural Sciences and
Engineering Research Council of Canada, and by the special Dean of
Science Fund at the University of Western Ontario.}
\thanks{$^\dag$Supported by the Mathematical Sciences Research
Institute, Berkeley}
\email{minac@uwo.ca}

\author[John Swallow]{John Swallow$^\ddag$}
\address{Department of Mathematics, Davidson College, Box 7046,
Davidson, North Carolina \ 28035-7046 \ USA}
\thanks{$^\ddag$Research supported in part by National Security
Agency grant MDA904-02-1-0061.} \email{joswallow@davidson.edu}

\begin{abstract}
In our previous paper we describe the Galois module structures of
$p$th-power class groups $K^\times/{K^{\times p}}$, where $K/F$ is
a cyclic extension of degree $p$ over a field $F$ containing a
primitive $p$th root of unity. Our description relies upon
arithmetic invariants associated with $K/F$. Here we construct
field extensions $K/F$ with prescribed arithmetic invariants, thus
completing our classification of Galois modules
$K^{\times}/K^{\times p}$.

\end{abstract}

\date\today

\maketitle

Let $F$ be a field of characteristic not $p$ containing a
primitive $p$th root of unity $\xi_p$.  For a cyclic field
extension $K/F$ with Galois group $\Gal(K/F)$ of order $p$, let $J
= {K^\times}/{K^{\times p}}$, and let $N$ denote the norm map from
$K$ to $F$.

In \cite{MS}, we proved that the structure of the
$\Fp[\Gal(K/F)]$-module $J$ is determined by the
following three arithmetic invariants:

\begin{itemize}
    \item $d=d(K/F):=\dim_{\Fp} {F^\times}/{N(K^\times)}$,
    \item $e=e(K/F):=\dim_{\Fp} N(K^{\times})/F^{\times p}$, and
    \item $\Upsilon(K/F):=$ $1$ or $0$ according to whether
    $\xi_p\in N(K^\times)$ or not.
\end{itemize}

Now if $G=\Z/p\Z$, then $J$ may be considered an $\Fp[G]$-module
via any isomorphism $G\cong \Gal(K/F)$, and the module structure
of $J$ is independent of the choice of isomorphism. It is a
fundamental problem to classify the isomorphism classes of modules
$J$ for all $K/F$ in our context. This problem is solved in
Theorem~\ref{th:mainth} below. Corollaries~1 and 2 of
Theorem~\ref{th:mainth} describe all modules $J$ in an explicit
way.

In the following theorem we determine the sets of invariants $(d,
e, \Upsilon)$ which may be realized by an extension $K/F$ and in
so doing classify all $\Fp[G]$-modules $J=J(K/F)$ up to
isomorphism.

\begin{theorem}\label{th:mainth}
    Let $p$ be a prime number. For arbitrary cardinal numbers $d$,
    $e$, and for $\Upsilon\in \{0,1\}$, there exists a cyclic
    field extension $K/F$ of degree $p$ containing a primitive
    $p$th root of unity with invariants $(d, e, \Upsilon)$ if and
    only if
    \begin{itemize}
        \item if $\Upsilon=0$, then $1\le d$,
        \item if $p>2$ then $1\leq e$, and
    \item if $p=2$ and $\Upsilon=1$ then $1\leq e$.
    \end{itemize}
\end{theorem}

From the theorem above and from \cite[Theorem~3 and
Corollary~2]{MS} we immediately obtain the following corollaries.
We denote by $M_{i,j}$ the $j$th cyclic module $\Fp[G]$ such that
$\dim_{\Fp}M_{i,j}=i$, where $j$ is a suitable index.

\begin{corollary} \label{co:maincor1}
    Let $p>2$ be a prime number, and let $G$ be a cyclic
    group of order $p$. Then an $\Fp[G]$-module $J$ is realizable
    as an $\Fp[\Gal(K/F)]$-module $K^{\times}/K^{\times p}$ for
    some cyclic $G$-extension $K/F$ such that $F$ contains a
    primitive $p$th root of unity if and only if there exist
    cardinal numbers $d$, $e$, and $\Upsilon\in\{0,1\}$ such that
    \begin{itemize}
        \item[$(i)$] If $\Upsilon=0$, then $1\le d$;
        \item[$(ii)$] $1\leq e$; and
        \item[$(iii)$]
        \begin{equation*}
            J=\left(\bigoplus_{j\in\Kf_1}M_{1,j}\right)
            \bigoplus\left(\bigoplus_{j\in\Kf_2}M_{2,j}\right)
            \bigoplus\left(\bigoplus_{j\in\Kf_p}
            M_{p,j}\right),
        \end{equation*}
        where
        \medskip
        \begin{enumerate}
            \item[$(1)$] $\vert\Kf_1\vert+1=2\Upsilon+d$,
            \item[$(2)$] $\vert\Kf_2\vert=1-\Upsilon$, and
            \item[$(p)$] $\vert\Kf_p\vert+1=e$.
        \end{enumerate}
    \end{itemize}
    The invariants $d$, $e$, and $\Upsilon$ determine the module $J$
    uniquely.
\end{corollary}

For $p=2$ using \cite[Theorem~3 and Corollary~3]{MS}, along with
our theorem above, we obtain the next corollary.

\begin{corollary} \label{co:maincor2}
    Now let $G$ be a cyclic group of order $2$. Then an
    $\F_2[G]$-module $J$ is realizable as an
    $\F_2[\Gal(K/F)]$-module $K^{\times}/K^{\times 2}$ for some
    quadratic extension $K/F$ with its arithmetic invariants
    $d(K/F)$, $e(K/F)$, and $\Upsilon(K/F)$ coinciding with $d$,
    $e$, and $\Upsilon\in\{0,1\}$, respectively, if and only if $d$, $e$, 
    and $\Upsilon$ satisfy the conditions below.

    \begin{itemize}
        \item if $\Upsilon=0$, then $1\leq d$,
    \item if $\Upsilon=1$, then $1\leq e$.
    \end{itemize}

    In this case
    \begin{equation*}
        J=\left(\bigoplus_{j\in\Kf_1}M_{1,j}\right)
        \bigoplus \left(\bigoplus_{j\in\Kf_2}M_{2,j}\right)
    \end{equation*}
    where
    \begin{enumerate}
        \item[(1)] $\vert\Kf_1\vert+1=2\Upsilon+d$ and
        \item[(2)] $\vert\Kf_2\vert+\Upsilon=e$.
    \end{enumerate}

    Moreover, such a module $J$ is determined uniquely by the
    invariants $2\Upsilon+d$ and by $e-\Upsilon$ if $e$ is finite
    and by $e$ alone if $e$ is infinite.
\end{corollary}

If $p>2$ then two $\Fp[G]$-modules are isomorphic if and only if
their invariants $d$, $e$, and $\Upsilon$ are the same, by
\cite[Corollary~2]{MS}. Thus we see in particular that if $p>2$,
then the arithmetic invariants of $J$ depend only upon the
isomorphism type of the $\Fp[G]$-module $J$.

In the case $p=2$, we see from \cite[Corollary~2]{MS} again that
the arithmetic invariants $d$, $e$, and $\Upsilon$ determine our
module $\F_2[G]$, but two isomorphic $\F_2[G]$-modules may have
different arithmetic invariants; see \cite[Corollary~3]{MS}. Here
is a very simple, concrete example illustrating this possibility.
Let $K_1/F_1$ be a quadratic extension of finite fields of
characteristic not $2$, $F_2=\R((t))$ be a field of power series
with coefficients in real numbers $\R$, and $K_2=F_2(\sqrt{-1})$.
Then both modules $K_{1}^{\times}/K_{1}^{\times 2}$ and
$K_{2}^{\times}/K_{2}^{\times 2}$ are isomorphic to a trivial
$\F_2[G]$-module $\F_2$, but their arithmetic invariants $(d_i,
e_i, \Upsilon_i)$ are $(0, 1, 1)$ for $i=1$ and $(2, 0, 0)$ for
$i=2$.

\section{Notation and Strategy} \label{S1}

In all that follows $F$ denotes a field, $F^{\times}=
F\setminus\{0\}$ the multiplicative group of $F$, $p$ a prime
number, and $F^{\times}/F^{\times p}$ the group of $p$th-power
classes of $F$. For each $f\in F^{\times}$ we denote by $[f]$ the
class of $f$ in $F^{\times}/F^{\times p}$. For each subset $A$ of
$F^{\times}$ we denote by $[A]$ the set of classes $\{[a]\ \vert\
a\in A\}$ and by $\mo[A]\mc$ the subgroup of $F^{\times}/F^{\times
p}$ generated by $[A]$.

We denote by $\xi_p$ a primitive $p$th root of unity in $F$. (Some
fields will be assumed to contain such a primitive $p$th root; for
the other fields in this paper, we will prove that a primitive
$p$th root is contained in the field.)  Observe that our
assumption that there exists a primitive $p$th root of unity
implies that $\chr(F)\neq p$.

For a Galois extension $K/F$, $\Gal(K/F)$ denotes the Galois group
and $N_{K/F}$ denotes the norm map from $K$ to $F$.  We denote by
$F^s$ the separable closure of $F$ and $G_F$ the absolute Galois
group $\Gal(F^s/F)$.  As usual, $H^i(G_F, \Fp)$ are Galois
cohomology groups of $F$ with coefficients in $\Fp$.  Since all
absolute Galois groups will be pro-$p$-groups, all considered
$\Fp$ modules are trivial.  Finally, let $\vert B \vert$ be the
cardinal number of a set $B$.

First observe that the conditions on $d$, $e$, and $\Upsilon$
listed in our theorem above are necessary:
\begin{enumerate}
    \item If $\Upsilon=0$, then $\xi_p\notin N_{K/F}(K^{\times})$
    and hence
    \begin{equation*}
        d =\dim_{\Fp} F^{\times}/N_{K/F}(K^{\times}) \geq 1;
    \end{equation*}
    \item If $p>2$ and $K=F(\root{p}\of{a})$ for a suitable $a\in
    F^{\times}\setminus F^{\times p}$, then $a=
    N_{K/F}(\root{p}\of{a})$ and hence
    \begin{equation*}
        e = \dim_{\Fp}N_{K/F}(K^{\times}) > 0.
    \end{equation*}
    \item If $p=2$, $\Upsilon=1$, and $K=F(\sqrt{a})$ for a suitable
    $a\in F^{\times}\setminus F^{\times 2}$, then $-1\in N_{K/F}(K^{\times})$
    since $\Upsilon=1$.  Consequently $a\in N_{K/F}(K^{\times})$,
    and thus $1\leq e$.
\end{enumerate}
Therefore in order to prove Theorem~\ref{th:mainth} when $p>2$,
it is sufficient to show, for each cardinal numbers $d$,
$e$ as above, for each $\Upsilon\in\{0,1\}$, and for each prime
number $p>2$, the existence of a field $F$ such that:
\begin{itemize}
    \item $F$ contains a primitive $p$th root $\xi_p$;
    \item $F^\times/F^{\times p}$ decomposes intro a direct sum of
    subgroups
    \begin{equation*}
        F^\times/F^{\times p} = D\oplus\mo[a]\mc\oplus E,
    \end{equation*}
    where $\dim_{\Fp}(\mo[a]\mc \oplus E) = e$ and, setting
    $K=F(\root{p}\of{a})$,
    \begin{enumerate}
        \item $[N_{K/F}(K^{\times})]=\mo[a]\mc \oplus E$;
        \item $\dim_{\Fp}(F^{\times}/N_{K/F}(K^{\times})) =
        \dim_{\Fp} D = d$; and
        \item $\Upsilon=0$ if and only if $\xi_p\notin
        N_{K/F}(K^{\times})$.
    \end{enumerate}
\end{itemize}

In the case $p=2$ and $\Upsilon=1$ we use the same conditions as
above, and if $p=2$ and $\Upsilon=0$ we require instead that $e =
\dim_{\F_2} E$ and $d=\dim_{\F_2}(D\oplus\mo[a]\mc)$. The latter
condition is imposed because $-1\notin N_{F(\sqrt{a})/F}
(F(\sqrt{a})^{\times})$ if and only if $a\notin N_{F(\sqrt{a})/F}
(F(\sqrt{a})^{\times})$.

Our strategy is to interpret the required conditions on
$F^{\times}/ F^{\times p}$ above in terms of Galois cohomology. We
then observe that these conditions are satisfied if $G_F$ is a
free product, in the category of pro-$p$-groups, of suitable
pro-$p$-groups $G_1$ and $G_2$, and finally we use the very nice
theorem proved by Efrat and Haran which guarantees the existence
of a field with $G_F$ above. This is one of the key results used
in our paper.

\begin{theorem} \label{th:efratharan}
    (Efrat-Haran; see \cite[Proposition~1.3]{EH}) Let
    $F_1,\dots,F_n$ be fields of equal characteristic such that
    $G_{F_1},\dots,G_{F_n}$ are pro-$p$-groups. Then there exists
    a field $F$ of the same characteristic such that
    \begin{equation*}
        G_F\cong G_{F_1}\star\cdots\star G_{F_n},
    \end{equation*}
    where the product is free in the category of pro-$p$-groups.
\end{theorem}

\noindent In order to apply this theorem, we show the existence of
the fields $F_1$ and $F_2$ such that $G_{F_1}$ and $G_{F_2}$ are
prescribed Galois groups $G_1$ and $G_2$.  We use the techniques
of henselian valuations and formal power series to construct
fields $F_1$ and $F_2$.

\section{Lemmas}

\subsection{Valued fields $F_1$ with prescribed residue field
$F_0$ and valuation group $\Gamma$}\

Let $v$ be a valuation on a field $F_1$, written additively. Then
we denote by $A_v$ the valuation ring $\{f\in F_1\ \vert\ v(f)\geq
0\}$; by $M_v$ the unique maximal ideal $\{f\in A_v\ \vert\
v(f)>0\}$ of $A_v$; by $F_v$ the residue field $A_v/M_v$ of $v$;
by $\Gamma$ the valuation group $v(F_1^{\times})$ of $v$; and by
$U$ the group $A_v\setminus M_v$ of units of $v$.

The following lemma is well known and we shall omit its
straightforward proof.

\begin{lemma} \label{le:pthcount}
    Let $F_1$ be a valued field with valuation $v$, valuation
    group $\Gamma=v(F_1^\times)$, and group of units $U$.
    For each prime $p\neq \chr(F_1)$ there exists an isomorphism
    \begin{equation*}
        \varphi:F_1^{\times}/F_1^{\times p}\longrightarrow U/U^{p}
        \oplus \Gamma/p\Gamma.
    \end{equation*}
    In particular
    \begin{equation*}
        \dim_{\Fp} F_1^{\times}/F_1^{\times p} = \dim_{\Fp}
        U/U^{p} + \dim_{\Fp}\Gamma/p\Gamma.
    \end{equation*}
\end{lemma}

It is well-known that for each field $F_0$ and for each
totally ordered abelian group $\Gamma$, there exists a field $F_1$
with a valuation $v:F_1\to\Gamma\cup\{\infty\}$ such that the
residue field $F_v$ is isomorphic to $F_0$ and the valuation group
is $\Gamma$.

In order to construct such a field, set
\begin{equation*}
    F_1 = F_0((\Gamma)) :=\{f:\Gamma\to F_0\ \vert\ \supp
    (f)\text{ is well-ordered}\}.
\end{equation*}
Thus a typical element $f\in F_1$ can be written as a formal sum
$f=\sum_{g\in\Gamma} a_g t^g$ such that the set $\supp (f):=\{g\in
G\ \vert\ a_g\neq 0\}$ is a well-ordered subset of $G$. The
valuation $v$ on $f$ is defined as: $v(0)=\infty$ and $v(f)= \min
\supp (f)$ for $f\neq 0$. An important property of the valued
field $F_1$ as above is the fact that it is henselian.
(See for example \cite[(1.3)]{Rib}.)
In what follows we will identify $F_v$ with $F_0$.
We will also assume that char $F_0\neq p$.

We will be particularly interested in controlling the $p$th-power
classes of such a field.  To do so, we choose particular groups
$\Gamma$ for our valuation groups.  These groups will be direct
sums of
\begin{equation*}
    \Z_{(p)}:=\left\{\frac{a}{b}\in\Q\ \big\vert\  a,b\in\Z,b \neq
    0; \ \text{if } a\neq 0 \text{ then } (a,b)=1, p\nmid
    b\right\}.
\end{equation*}
Observe that $\Z_{(p)}$ is the valuation ring of a $p$-adic
valuation on $\Q$. Let $I$ be any non-empty, well-ordered set.
Then set
\begin{equation*}
    \Gamma=\Z_{(p)}^{(I)} := \left\{\gamma: I\to \Z_{(p)}\ \Big|
    \ \vert\supp(\gamma)\vert<\infty\right\}.
\end{equation*}
Thus $\Gamma$ is a direct sum of $\vert I\vert$ copies of
$\Z_{(p)}$. Observe that $\Z_{(p)}$ carries a natural ordering
induced from $\Q$, and then we may order $\Gamma$
lexicographically, as follows. Let
$\gamma_1\neq\gamma_2\in\Gamma$. Then $\gamma_1<\gamma_2$ if and
only if $\gamma_1(i)<\gamma_2(i)$ for the least element $i\in I$
such that $\gamma_1(i)\neq\gamma_2(i)$. Then $\Gamma$ is a
linearly ordered abelian group. Recall that each non-empty set can
be well-ordered. (See \cite[Appendix~2, Theorem~4.1]{La}.)

We choose $\Gamma$ as above because $G_{F_0((\Gamma))}$ will be
pro-$p$ (see Lemma~\ref{le:gf1prop} below) and because we may
control the $p$th-power classes with the following lemma. This
well-known lemma follows from Lemma~\ref{le:pthcount} and the fact
that the valued field $F_1$ is henselian. It is also an immediate
consequence of \cite[Lemma~1.4]{W}. Therefore we shall omit its
proof.

\begin{lemma}\label{le:pthcountf1}
    Let $F_1=F_0((\Gamma))$ as above.  Then
    \begin{align*}
        \dim_{\Fp} F_1^\times/F_1^{\times p} &= \dim_{\Fp}
        F_0^\times/F_0^{\times p} + \dim_{\Fp} \Gamma/p\Gamma
        \\ &= \dim_{\Fp} F_0^\times/F_0^{\times p} +
        \vert I\vert.
    \end{align*}
\end{lemma}

Finally, we record a criterion for $G_{F_1}$ being pro-$p$:

\begin{lemma} \label{le:gf1prop}
    Let $F_1=F_0((\Gamma))$ as above with $\chr(F_0) = 0$ and
    $G_{F_0}$ pro-$p$.  Then $G_{F_1}$ is pro-$p$ as well.
\end{lemma}

\begin{proof}
    From basic valuation theory, nicely summarized in \cite[pages~3
    and 4]{K}, and the fact that $F_1$ above is henselian,
    \begin{equation*}
        G_{F_1}\cong T\rtimes G_{F_0},
    \end{equation*}
    where the action of $G_{F_0}$ on $T$ is uniquely determined by
    the cyclotomic character mapping $G_{F_0}$ into a group of
    automorphisms of a group of roots of unity contained in
    $F_0^s$, and $T\cong\Z_{p}^{I}$, the topological
    product of $\vert I\vert$ copies of $\Z_p$. In particular, if
    $G_{F_0}$ is a pro-$p$-group, so is $G_{F_1}$.
\end{proof}

\subsection{$H^2(G_{F_1}, \Fp)$ for henselian valued fields $F_1$}\

Now we study $H^2(G_{F_1}, \Fp)$ for our henselian valued fields
$F_1$.  The next lemma, taken from \cite{W}, will be used in the
proof of Theorem~\ref{th:mainth} to show that norm groups of
cyclic $p$-extensions of $F_1$ are not too large.

Suppose that $F_1$ is a field endowed with a henselian valuation
$v$ with valuation group $\Gamma=v(F_1^\times)$. Let $F_1^{{nr}}$
denote the maximal unramified extension of $F_1$ in its separable
closure $F_1^s$. Then $G_{F_0}\cong
G_{F_1}/\Gal(F_1^s/F_1^{{nr}})$.  Therefore, after identifying
these groups, we have the inflation map
\begin{equation*}
    \inf = \inf{}_{F_0}^{F_1} : H^{*}(G_{F_0}, \Fp)
    \longrightarrow H^{*}(G_{F_1}, \Fp).
\end{equation*}
(See \cite[page~483]{W}.)

Moreover, from basic Kummer theory we have the canonical
isomorphism
\begin{equation*}
    \varphi_F : F^{\times}/F^{\times p}\longrightarrow
    H^{1}(G_F, \Fp),
\end{equation*}
as well as the corresponding canonical isomorphisms
$\varphi_{F_i}$, $i=1$, $2$.  We will denote by $(f)_F$ or
$(f_i)_{F_i}$ the images $\varphi_F([f])$ or
$\varphi_{F_i}([f_i])$.  If the context is clear we will omit the
subscript.

Assume next that $\{\pi_j,\ j\in \Jc\}$ is a set of elements of
$F_1^{\times}$ such that their images in $\Gamma/p\Gamma$ form a
basis of $\Gamma/p\Gamma$ over $\Fp$. Then we have the following
lemma, obtained as a special case of a theorem of Wadsworth.

\begin{lemma} \label{le:h2ofhens} \cite[Theorem~3.6, page~483]{W}.
    Let $F_1$ and $F_0$ be as above. Then
    \begin{align*}
        H^2(G_{F_1}, \Fp) = &\inf(H^2(G_{F_0},\Fp))\oplus_{j\in \Jc}
        \left( \inf(H^1(G_{F_0}, \Fp))\cup(\pi_j)\right) \\
        &\oplus_{\{j_1,j_2\}\subset \Jc,\ j_1\neq j_2}
        \left((\pi_{j_1})\cup(\pi_{j_2})\right).
    \end{align*}
    Moreover, for each $j\in \Jc$,
    \begin{equation*}
        \inf(H^1(G_{F_0},\Fp)) \cong \inf(H^1(G_{F_0},\Fp))
        \cup(\pi_j)
    \end{equation*}
    and for each $j_1,j_2\in \Jc$ such that $j_1\neq
    j_2$, we have $(\pi_{j_1})\cup(\pi_{j_2})\neq 0$.
\end{lemma}

Note that in the last summand of Lemma~\ref{le:h2ofhens} the sum
ranges over subsets $\{j_1,j_2\},j_1\neq j_2$ of $\Jc$ and a
choice between $(\pi_{j_{1}})\cup(\pi_{j_{2}})$ and
$(\pi_{j_{2}})\cup(\pi_{j_{1}})$ is arbitrary but fixed.

\subsection{Residue fields $F_0$ with prescribed absolute
Galois group}\

In our construction of $F$ we choose a residue field $F_0$
depending on $\Upsilon$ and $p$.  If $\Upsilon=1$ we will simply
put $F_0=\C$, but when $\Upsilon=0$ we require some special
properties of $F_0$. In particular, in order that our cyclic
extension $K=F(\root{p}\of{a})$ have the desired invariant
$\Upsilon=0$, we require that $K$ does not embed in a cyclic
Galois extension $L$ over $F$ with degree $[L:F]=p^2$, for this is
equivalent to $\xi_p\notin N_{K/F}(K^\times)$ by
\cite[Theorem~3]{A}.

To ensure that this nonembeddability condition holds, as well as
to ensure that a certain nonabelian group of order $p^3$ does not
occur as a Galois group over the field, we choose residue fields
$F_0$ with absolute Galois groups taking a special form, and it is
also convenient to require that $\vert F_0^\times/F_0^{\times
p}\vert$ is small. As it turns out, we may choose some suitable
algebraic infinite extension of $\Q$. Finitely generated
pro-$p$-absolute Galois groups over $\Q$ and more generally any
global field, were nicely classified in \cite{E2}. (See also
\cite{E1} and \cite{JP} for related results and techniques.)

The extensions we will need for $p>2$ are given in the following

\begin{lemma} \label{le:specextsofQ} \cite[page~84]{E2}
    For each prime $p>2$ there exists an algebraic extension
    $F_{0,p}$ of $\Q$ such that
    \begin{equation*}
        G_{F_{0,p}} = \mo\sigma,\tau\ \vert\
        \sigma\tau\sigma^{-1}= \tau^{p+1}\mc_{\text{pro-$p$}}
    \end{equation*}
    where the presentation is in the category of pro-$p$-groups.
\end{lemma}

Observe that the maximal abelian extension $F_{0,p}^{ab}$ of
$F_{0,p}$ has $G_{F_{0,p}}^{ab} := \Gal(F_{0,p}^{ab}/ F_{0,p})$
equal to
\begin{equation*}\label{eq:gfopab}
    G_{F_{0,p}}^{ab} \cong \mo\bar\sigma,\bar\tau\ \vert
    \ \bar\tau^{p}=1\mc=\Z_p\times\Z/p\Z
\end{equation*}
for all primes $p>2$.

\subsection{Field arithmetic and free pro-$p$ products}\

In this section we collect lemmas giving information about a field
$F$ derived from the structure of $G_F$, especially when $G_F$ is
a free pro-$p$ product of two pro-$p$ groups $G_{F_1}$ and
$G_{F_2}$.

First we record a lemma detecting the presence of primitive $p$th
roots of unity in a field $F$, based only on the structure of
$G_F$.

\begin{lemma} \label{le:pthrootexists}
    Suppose that $p>2$ and that $F$ is a field with $\chr(F)\neq
    p$ and $G_F$ pro-$p$. Then $\xi_p\in F^{\times}$.
\end{lemma}

\begin{proof}
    Because $\chr(F)\neq p$, there exists a primitive $p$th root
    $\xi_p$ of unity in $F^s$. If $\xi_p\in F^s\setminus F$ then
    $F(\xi_p)/F$ is a nontrivial Galois extension of degree
    $[F(\xi_p):F]<p$. Therefore $G_F$ has a nontrivial finite
    quotient of order coprime with $p$. This contradicts our
    assumption that $G_F$ is a pro-$p$-group. Hence $\xi_p\in
    F^{\times}$ as asserted.
\end{proof}

Now suppose that $G_F=G_{F_1}\star G_{F_2}$ for pro-$p$ absolute
Galois groups $G_F$, $G_{F_1}$, and $G_{F_2}$, where the free
product is taken in the category of pro-$p$-groups.  From
\cite[(4.3) Satz]{N} we see that the restriction homomorphism
\begin{equation}\label{eq:res}
    \res : H^{1}(G_F, \Fp)\longrightarrow H^{1}(G_{F_1},
    \Fp)\oplus H^{1}(G_{F_2}, \Fp)
\end{equation}
is an isomorphism.  Now given $(f)_F$ in $H^1(G_F, \Fp)$, we
denote the image $\res (f)_F$ by
\begin{equation*}
    \res (f)_F = (f)_{G_{F_1}} \oplus (f)_{G_{F_2}}.
\end{equation*}
This notation distinguishes, then, between $(f)_{F_1}$, which
denotes $\varphi_{F_1}(f)$ for $f\in F_1^\times$, and
$(f)_{G_{F_1}}$, which denotes the projection of $\res
\varphi_{F}(f)$ onto the first summand.

One way of interpreting this restriction map is with the following

\begin{lemma}\label{le:resinterp}
    Let $G_F=G_{F_1}\star G_{F_2}$ be pro-$p$ absolute Galois
    groups of fields containing a $p$th root of unity, and
    suppose that we have the following sequence:
    \[
    \xymatrix{ G_F\ar@{>>}[r]^{\can\ } & G_{F_1} \ar@{>>}[r] &
    \Z/p\Z. }
    \]
    Here the canonical map $\can$ is an identity on $G_{F_1}$ and
    contains $G_{F_2}$ in its kernel.

    Then the right-hand surjection and the composed surjection
    correspond to fields $K_1=F_1(\root{p}\of{a_1})$ and
    $K=F(\root{p}\of{a})$, respectively, where
    $(a)_{G_{F_1}}=(a_1)_{F_1}$ and $(a)_{G_{F_2}}=0$.
\end{lemma}

\begin{proof}
    The surjections are continuous homomorphisms, hence elements
    of $H^1(G_{F_1},\Fp)$ and $H^1(G_{F},\Fp)$, respectively, and
    the right-hand surjection is clearly the restriction of the
    composed surjection.  The remainder follows by Kummer theory.
\end{proof}

\begin{lemma} \label{le:cupprod}
    Let $G_F=G_{F_1}\star G_{F_2}$ be pro-$p$ absolute Galois
    groups and suppose that $(a)_F$ and $(b)_F$ satisfy
    $(a)_{G_{F_1}} = (b)_{G_{F_2}} = 0$.  Then
    \begin{equation*}
        (a)_F \cup (b)_F = 0 \in H^{2}(G_F, \Fp).
    \end{equation*}
\end{lemma}

The lemma follows from \cite[(4.1) Satz]{N} and from
\cite[Prop.~7.3, page~191]{Ris}. However, we prove our lemma by
translating the cup products into obstructions to basic Galois
embedding problems, yielding an interesting Galois-theoretic
variant of the proof.

\begin{proof}
    If $(a)=0$ or $(b)=0$, we are done.  Otherwise, the conditions
    $(a)_{G_{F_1}}=(b)_{G_{F_2}}=0$ imply that $(a)$ and $(b)$ are
    linearly independent in $H^1(G_F, \Fp)$.

    Now let $H_{p^3}$ be the Heisenberg group of
    order $p^3$:
    \begin{align*}
        H_{p^3} = \mo &v_1, v_2, w\ \vert\ v_1^p=v_2^p=w^p=1,
        v_2v_1=wv_1v_2, \\ &[v_1,w]=[v_2,w]=1 \mc
    \end{align*}
    In the case $p=2$, $H_{8}$ is the familiar dihedral group
    $D_4$.

    By \cite[Corollary, page 523 and Theorem 3(A)]{M}, if $(a)$
    and $(b)$ are linearly independent, then $(a)\cup (b)=0$ if
    and only if $H_{p^3}$ is the Galois group $\Gal(M/F)$ of a
    Galois extension $M$ of $F$ containing $F(\root{p}\of{a},
    \root{p}\of{b})$ in such a way that
    \begin{equation*}
        H_{p^3}/\mo v_1,w\mc = \Gal(F(\root{p}\of{a})/F)
        \text{ and }
        H_{p^3}/\mo v_2,w\mc = \Gal(F(\root{p}\of{b})/F).
    \end{equation*}

    Now consider the commutative diagram
    \[
    \xymatrix{ & & G_F\ar@{>>}[ld]_{\delta_1} \ar@{.>>}[dd]^\beta
    \ar@{>>}[rd]^{\delta_2} & & \\ & G_{F_1} \ar@{>>}[dl]_{\alpha_1}
    & & G_{F_2} \ar@{>>}[dr]^{\alpha_2} & \\ \Z/p\Z
    \ar@{^{(}->}[rr]^{1\mapsto v_1} & & H_{p^3}
    \ar@(dr,dl)[rr]_{\mo v_1, w\mc\mapsto 0; \ v_2\mapsto 1}
    \ar@(dl,dr)[ll]^{\mo v_2, w\mc\mapsto 0; \ v_1\mapsto 1} & & \Z/p\Z
    \ar@{_{(}->}[ll]_{1\mapsto v_2}
    }
    \]
    Let $a_2\in F_2^\times$ and $b_1\in F_1^\times$ satisfy
    $(a_2)_{F_2} = (a)_{G_{F_2}}$ and $(b_1)_{F_1} =
    (b)_{G_{F_1}}$.  Then set $K_1=F_1(\root{p}\of{b_1})$ and
    $K_2=F_2(\root{p}\of{a_2})$; these are $\Z/p\Z$-extensions of
    $F_1$ and $F_2$, respectively.  We may then identify the
    left-hand $\Z/p\Z$ in the diagram with $\Gal(K_1/F_1)$ so that
    $\alpha_1$ is the surjection of Galois theory.  Similarly, the
    right-hand $\Z/p\Z$ may be identified with $\Gal(K_2/F_2)$
    with $\alpha_2$ the surjection of Galois theory.  Finally,
    the topmost surjections $\delta_1$ and $\delta_2$ are canonical.

    By Lemma~\ref{le:resinterp}, the surjections
    $\delta_1\alpha_1$ and $\delta_2\alpha_2$ correspond to fields
    $F(\root{p}\of{b})$ and $F(\root{p}\of{a})$, respectively.
    Now because $G_F=G_{F_1}\star G_{F_2}$, there exists a
    homomorphism $\beta: G_F\to H_{p^3}$, and because $v_1$ and
    $v_2$ generate $H_{p^3}$, $\beta$ is a surjection.

    Hence $H_{p^3}$ is a Galois group over $F$ corresponding to a
    normal subgroup $H$ of $G_F$.  Consider the smallest normal
    subgroup $H_1$ of $G_{F}$ containing $H$ and $G_{F_2}$.  Then
    by the diagram, $G_F/H_1$ is the left-hand $\Z/p\Z$, which
    corresponds to $F(\root{p}\of{b})$, and $H_1/H$ is $\mo v_2,
    w\mc$.  Now consider the smallest normal subgroup $H_2$ of
    $G_{F}$ containing $H$ and $G_{F_1}$.  Then by the diagram,
    $G_F/H_2$ is the right-hand $\Z/p\Z$, which corresponds to
    $F(\root{p}\of{a})$, and $H_2/H$ is $\mo v_1, w\mc$.

    Hence $(a)\cup (b)=0$.
\end{proof}

Finally, we close with with a companion to
Lemma~\ref{le:resinterp}. In Lemma~\ref{le:abeloffree} below, $\pi$
denotes the canonical homomorphism of $G$ to $G_1$ which is an identity
on $G_1$, and is trivial on $G_2$.

\begin{lemma} \label{le:abeloffree}
    Let $G=G_1\star G_2$ be a free product of $G_1$ and $G_2$ in
    the category of pro-$p$-groups. Suppose that $A\cong \Z/p\Z$
    is a factor group of $G_1$ such that the surjection $G_1\to
    \Z/p\Z$ does not factor through $\Z/p^2\Z$. Then the following
    commutative diagram cannot occur:
    \[
    \xymatrix{ G\ar@{>>}[r]^-{\pi} \ar@{>>}[d] & G_1 \ar@{>>}[d] \\
    \Z/p^2\Z \ar@{>>}[r] & A. }
    \]
\end{lemma}

\begin{proof}
    Suppose that contrary to our statement, such a diagram as the
    above exists. Then by passing to quotients by commutator
    subgroups we obtain
    \[
    \xymatrix{ G^{ab} \ar@{>>}[r]^-{\pi^{ab}}
    \ar@{>>}[d]_{\alpha} &
    G_1^{ab} \ar@{>>}[d]^{\gamma} \\
    \Z/p^2\Z \ar@{>>}[r]^{\beta} & A. }
    \]
    But $G^{ab}\cong G_1^{ab}\times
    G_2^{ab}$ and the canonical surjection onto
    $G_1^{ab}$ is given by the projection map.  Let
    $\delta$ be a splitting map of the projection map.  Then
    $\gamma = \beta\alpha\delta$, contradicting the hypothesis.
\end{proof}

\section{Proof of the Theorem} \label{se:theproof}

First we define fields $F_0$, $F_1$, $F_2$, and $F$ using our
given cardinal numbers $d$, $e$, and $\Upsilon$, as well as the
prime number $p$.  Then we define the cyclic Galois extension
$K/F$ of degree $p$ and check that the arithmetic invariants of
$K/F$ coincide with $d$, $e$, and $\Upsilon$.

\subsection{Constructing $F_0$, $F_1$, $F_2$, and $F$}\

If $\Upsilon=1$ then let $F_0=\C$. If $\Upsilon=0$ and $p=2$, let
$F_0=\R$. (See Proposition~\ref{pr:trandeg} for alternative
choices in these two cases.)  If $\Upsilon=0$ and $p>2$ then let
$F_{0,p}$ be the algebraic extension of $\Q$ of
Lemma~\ref{le:specextsofQ}. In the first two cases we see
trivially that $\xi_p\in F_0$, and in the last case $\xi_p\in F_0$
by Lemma~\ref{le:pthrootexists}. Observe that
\begin{equation}\label{eq:f0pth}
        \dim_{\Fp} F_{0}^{\times}/F_{0}^{\times p}=
    \begin{cases}
        0, \text{ if } \Upsilon=1;\\
        1, \text{ if } \Upsilon=0, \ p=2;\\
        2, \text{ if } \Upsilon=0, \ p>2.
    \end{cases}
\end{equation}

We next construct the field $F_1$. Because $\Upsilon=0$ implies
$1\leq d$, for either choice of $\Upsilon\in \{0,1\}$ there exists
a well-ordered set $I_1$ such that $\vert I_1
\vert+1=d+2\Upsilon$. Let $\Gamma_1=\Z_{(p)}^{(I_1)}$ be a direct
sum of $\vert I_1 \vert$ copies of $\Z_{(p)}$. Then $\Gamma_1$ is
a linearly ordered abelian group. Finally set $F_1 :=
F_0((\Gamma_1))$. From Lemmas~\ref{le:pthcountf1} and
\ref{le:gf1prop} it follows that $G_{F_1}$ is a pro-$p$-group and
\begin{equation*}
    \dim_{\Fp}F_{1}^{\times}/F_{1}^{\times p}=
    \dim_{\Fp}F_{0}^{\times}/F_{0}^{\times p} + \vert I_1 \vert.
\end{equation*}
Hence
\begin{equation}\label{eq:f1pth}
    \dim_{\Fp} F_{1}^{\times}/F_{1}^{\times p}=
    \begin{cases}
    d+1, &\text{if } \Upsilon=1;\\
       d, &\text{if } \Upsilon=0, \ p=2;\\
       d+1, & \text{if } \Upsilon=0, \ p>2.
    \end{cases}
\end{equation}

Similarly, we construct $F_2$ as follows.  Because $p>2$ and also
$p=2$ and $\Upsilon=1$ implies
$e>0$, there exists a well-ordered set $I_2$ such that $1+\vert
I_2 \vert = e$ in either of the cases $p>2$ or $\Upsilon=1$,
$p=2$, and $\vert I_2 \vert=e$ in the case $\Upsilon=0$, $p=2$.
Then again $\Gamma_2=\Z_{(p)}^{(I_2)}$ is a linearly ordered
abelian group. We set $F_2 := \C((\Gamma_2))$. Then from
\cite[pages~3~and~4]{K} it follows that
$G_{F_2}\cong\Z_{p}^{I_2}$, the topological product of $\vert I_2
\vert$ copies of $\Z_p$. In particular
\begin{equation}\label{eq:edef}
        e =
    \begin{cases}
        \dim_{\Fp} F_{2}^{\times}/F_{2}^{\times p}+1, &\text{if }
        p>2 \text{ or } p=2 \text{ and } \Upsilon=1;\\
        \dim_{\Fp} F_{2}^{\times}/F_{2}^{\times p}, &
        \text{if }
        p=2, \ \Upsilon=0.\\
    \end{cases}
\end{equation}

From Theorem~\ref{th:efratharan} we see that there exists a field $F$ of 
characteristic zero, such
that $G_F=G_{F_1}\star G_{F_2}$ is a free product of $G_{F_1}$ and
$G_{F_2}$ in the category of pro-$p$-groups. In particular $G_F$
is again a pro-$p$-group, and from Lemma~\ref{le:pthrootexists} we
see that $F$ contains a primitive $p$th root of unity.

\subsection{Constructing $K/F$}\label{kf}\

We define the cyclic extension $K/F$ of degree $p$ as
$K=F(\root{p}\of{a})$ where $[a]\in F^{\times}/F^{\times p}$ is
chosen via the isomorphism \eqref{eq:res}.

If $\Upsilon=1$ then let $(a_1)_{F_1}$ be any nontrivial element
in $H^1(G_{F_1}, \Fp)$, which is possible since $1\leq \dim_{\Fp}
H^1(G_{F_1}, \Fp) = d+1$.  Let $(a_2)_{F_2}=0$ in $H^2(G_{F_2},
\Fp)$.

Now suppose $\Upsilon=0$ and $p>2$.  We denote the fixed field of
the factor $\Z_p$ in $G_{F_0}^{ab}\cong\Z_p\times\Z/p\Z$ acting on
a maximal abelian extension $F_{0}^{ab}$ of $F_0$ as
$K_0:=F_0(\root{p}\of{b})$. (See the discussion following
Lemma~\ref{le:specextsofQ}.) Because $F_{0}^{\times}/F_{0}^{\times
p}$ is naturally isomorphic to a subgroup of
$F_{1}^{\times}/F_{1}^{\times p}$, we may set $(a_1)_{F_1} = \inf
(b)_{F_0}\neq 0$ in $H^1(G_{F_1}, \Fp)$ and $(a_2)_{F_2}=0$ in
$H^2(G_{F_2},\Fp)$.

Finally assume that $\Upsilon=0$ and $p=2$.  We set
$(a_1)_{F_1}=(-1)_{F_1}\neq 0$ in $H^1(G_{F_1}, \F_2)$ and
$(a_2)_{F_2}=0$ in $H^2(G_{F_2},\F_2)$.

Now by \eqref{eq:res} there exists $a\in F^\times$ such that
$(a)_F\neq 0$ and $(a)_{G_{F_i}}=(a_i)_{F_i}$, $i=1$, $2$.  Since
we have already determined that $\xi_p\in F^{\times}$, we see that
$K=F(\root{p}\of{a})$ is a cyclic extension of degree $p$. Hence
it remains to show that the arithmetic invariants $d(K/F)$,
$e(K/F)$, and $\Upsilon(K/F)$ coincide with the prescribed
cardinal numbers $d$, $e$, and $\Upsilon$ respectively.

\subsection{Determining $d(K/F)$ and $e(K/F)$ via annihilators}\

We first observe some relationships among cup products of elements
in $H^1(G_F, \Fp)$.  Recall that Lemma~\ref{le:cupprod} tells us
that if $c_1, c_2\in F^\times$ with $(c_1)_{G_{F_1}} =
(c_2)_{G_{F_2}} = 0$ then
\begin{equation}\label{eq:f1f2}
    (c_1)_F \cup (c_2)_F = 0\in
    H^2(G_F, \Fp).
\end{equation}
Moreover, recall that by \cite[Satz~4.1]{N}
\begin{equation*}
    H^2(G_F,\Fp)\cong H^2(G_{F_1},\Fp)\times H^2(G_{F_2},\Fp),
\end{equation*}
where the isomorphism is induced by the restriction maps
\begin{equation*}
    \res_i : H^2(G_F,\Fp) \to H^2(G_{F_i}, \Fp), i=1, 2.
\end{equation*}
Because these restriction maps commute with cup product maps
(\cite[Proposition~7.3, page~191]{Ris}), we see that if $c_1,
c_2\in F^\times$ such that $(c_1)_{G_{F_2}}=(c_2)_{G_{F_2}}=0$ in
$H^2(G_{F_2}, \Fp)$ then
\begin{equation}\label{eq:tof1}
    (c_1)_F \cup (c_2)_F \neq 0\quad
    \text{iff}\quad (c_1)_{G_{F_1}}\cup (c_2)_{G_{F_1}}\neq 0,
\end{equation}
where the first cup product lies in $H^1(G_F, \Fp)$ and the second
cup product lies in $H^1(G_{F_1}, \Fp)$.  By symmetry the
analogous statement with $F_1$ and $F_2$ exchanged holds as well.

Now we calculate $e(K/F)$. We adopt the following notation for an
annihilator of $\phi\in H^1(G_{F},\Fp)$:
\begin{equation*}
    \ann_{F_1}\phi := \{\eta\in H^1(G_{F}, \Fp)\ \vert\
    \eta\cup\phi=0\}.
\end{equation*}
For $x\in F^{\times}$ we have $x\in
N_{K/F}(K^{\times})$ if and only if $(x)_F \cup (a)_F =0$ in
$H^2(G_F, \Fp)$.  Therefore, using \eqref{eq:f1f2} and
\eqref{eq:tof1} with the fact that in every case
$(a)_{G_{F_2}}=0$,
\begin{align*}
    e(K/F) &= \dim_{\Fp} N(K^{\times})/F^{\times p} \\
        &= \dim_{\Fp}\{(x)_F \in H^1(G_F,\Fp)\ \vert \ (a)_F
        \cup (x)_F=0\} \\
        &= \dim_{\Fp}\{(y)_{F_1} \in H^1(G_{F_1}, \Fp)\ \vert
        \ (a)_{G_{F_1}} \cup(y)_{F_1}=0\}\\
        &\phantom{=\ \ } + \dim_{\Fp}
        H^1(G_{F_2},\Fp) \\
        &= \dim_{\Fp} \ann_{F_1}(a)_{G_{F_1}} + \dim_{\Fp}
        F_2^\times/F_2^{\times p}.
\end{align*}

If $\Upsilon=1$ then $H^1(G_{F_0}, \Fp)=0$, and so by
Lemma~\ref{le:h2ofhens} we have $\dim_{\Fp} \ann_{F_1}
(a)_{G_{F_1}} = 1$.

If $p>2$ and $\Upsilon=0$ then by Lemma~\ref{le:h2ofhens},
\begin{equation*}
    \dim_{\Fp} \ann_{F_1}(a)_{G_{F_1}} = \dim_{\Fp} \ann_{F_0}
    (b)_{F_0},
\end{equation*}
where $b$ was chosen so that $K_0=F_0(\root{p}\of{b})$ was the
fixed field of the factor $\Z_p$ in $G_{F_0}^{ab}$. Now because
$(b)_{F_0} \cup (b)_{F_0} = 0$ as an identity in $H^2(G_F, \Fp)$
for $p>2$,
\begin{equation*}
    1\le \dim_{\Fp} \ann_{F_0} (b)_{F_0} \le 2.
\end{equation*}
On the other hand, let $c\in F^\times\setminus F^{\times p}$ be
defined so that $F_0(\root{p}\of{c})$ is contained in
the fixed field of the
factor $\Z/p\Z$ in $G_{F_0}^{ab}$.  Then by \eqref{eq:f0pth}, $\{
(b)_{F_0}, (c)_{F_0}\}$ spans $H^1(G_{F_0}, \Fp)$. Furthermore,
since by Lemma~\ref{le:specextsofQ}, $H_{p^3}$ is not a quotient
of $G_{F_{0,p}}$, $(b)_{F_0} \cup (c)_{F_0}\neq 0$ by
\cite[Corollary, page 523 and Theorem 3(A)]{M}. We conclude that
$\dim_{\Fp} \ann_{F_0} (b)_{F_0} = 1$.

Finally, if $p=2$ and $\Upsilon=0$ then again by
Lemma~\ref{le:h2ofhens}
\begin{equation*}
    \dim_{\F_2} \ann_{F_1}(a)_{G_{F_1}} = \dim_{\F_2}
    \ann_{F_0} (-1)_{F_0}.
\end{equation*}
As $F_0=\R$, $\dim_{\F_2} F_0^\times/F_0^{\times 2} = 1$ and
$(-1)_{\R}\cup (-1)_{\R} \neq 0$, yielding $\dim_{\F_2} \ann_{F_1}
(a)_{G_{F_1}} = 0.$

Combining our results with \eqref{eq:edef}, we have that
$e(K/F)=e$.

Now we turn to a similar calculation of $d(K/F)$.  Again using
\eqref{eq:f1f2} and \eqref{eq:tof1} with the fact that in every
case $(a)_{G_{F_2}}=0$,
\begin{align*}
    d(K/F)  &= \dim_{\Fp}(H^1(G_F,\Fp)/ \ann_{F} (a)_F)  \\
    &= \dim_{\Fp}(H^1(G_{F_1},\Fp)/\ann_{F_1} (a)_{G_{F_1}}) \\
    &= d.
\end{align*}
For this last equality we use \eqref{eq:f1pth} together with
calculations of the dimension of $\ann_{F_1} (a)_{G_{F_1}}$
already achieved. Hence $d(K/F)=d$ in all cases.

\subsection{Determining $\Upsilon$ via quotients of $G_F$}\

It remains to show that $\Upsilon(K/F)=\Upsilon$. First consider
the case $\Upsilon=1$.  Since $F_0=\C$, $\xi_p$ is a $p$th power
in $F_1$, and $F_1(\root{p}\of{a})$ embeds in a
$\Z/p^2\Z$-extension $F_1(\root{p^2}\of{a})$ of $F_1$. Then the
surjection $G_{F_1}\to\Gal(K_1/F_1)$ factors through $\Z/p^2\Z$.
Following the surjection with the canonical surjection $G_F\to
G_{F_1}$, we see that $\Z/p^2\Z$ is a factor group of $G_F$.
Moreover, by Lemma~\ref{le:resinterp}, the surjection $G_F \to
\Z/p\Z$ corresponds to $K$.  Hence $K/F$ embeds in a
$\Z/p^2\Z$-extension of $F$.  By \cite[Theorem~3]{A}, $\xi_p\in
N_{K/F}(K^{\times})$. Therefore $\Upsilon(K/F)=1=\Upsilon$.

Now consider the case $\Upsilon=0$ and $p>2$.  Because
$K_0=F_0(\root{p}\of{b})$ does not embed in a $\Z/p^2\Z$-extension
of $F_0$, $\xi_p\notin N_{K_0/F_0}(K_0^\times)$. (See \ref{kf} for
the definition of $K_0$.) Hence $(b)_{F_0} \cup(\xi_p)_{F_0}\neq
0$ in $H^1(G_{F_0},\Fp)$, and, by Lemma~\ref{le:h2ofhens}, $(\inf
(b))\cup (\xi_p)_{F_1} \neq 0$ in $H^1(G_{F_1},\Fp)$ as well.
Choose $b_1\in F_1^\times$ so that $(b_1)_{F_1} = \inf (b)_{F_0}$
and set $K_1=F_1(\root{p}\of{b_1})$. Then $\xi_p \notin
N_{K_1/F_1} (K_1^\times)$ and by \cite[Theorem~3]{A} the field
extension $K_1/F_1$ does not embed in a $\Z/p^2\Z$-extension of
$F_1$. Therefore the surjection $G_{F_1}\to
\Gal(F_1(\root{p}\of{b})/F_1)$ does not factor through $\Z/p^2\Z$.
From Lemmas~\ref{le:resinterp} and \ref{le:abeloffree} we see that
surjection $G_F\to\Gal(K/F)$ does not factor through $\Z/p^2\Z$.
Again by \cite[Theorem~3]{A} we conclude that $\xi_p\notin
N_{K/F}(K^{\times})$ and therefore $\Upsilon(K/F)=0=\Upsilon$.

Finally consider the case $\Upsilon=0$ and $p=2$. Because
$-1\notin N_{\C/\R}(\C)$ from Lemma~\ref{le:h2ofhens} we conclude
that $-1\notin N_{K_1/F_1}(K_{1}^{\times})$.  By
\cite[Theorem~3]{A} the surjection $G_{F_1}\to\Gal(K_1/F_1)$ does
not factor through $\Z/4\Z$. As in the previous case, by
Lemmas~\ref{le:resinterp} and \ref{le:abeloffree} we see that
$K/F$ does not embed in a $\Z/4\Z$-extension of $F$. Again by
\cite[Theorem~3]{A} we see that $-1\notin N_{K/F}(K^{\times})$ and
therefore $\Upsilon(K/F)=0=\Upsilon$.

Hence we have checked in all cases that $\Upsilon(K/F)=\Upsilon$,
and our proof of Theorem~\ref{th:mainth} is now complete. \qed

\begin{remark*}
In \cite[Lemma~1.2 and Proposition~1.3]{EH} Efrat and Haran
construct some fields with prescribed absolute Galois groups
together with some bounds on the transcendence degrees of these
fields. These bounds, together with the replacement of $\C$ by an
algebraic closure of $\Q$ and of $\R$ by a real-closed algebraic
number field $R$ in the case when $p=2$ and $\Upsilon=0$ in our
proof above, yield the following proposition.
\end{remark*}

\begin{proposition} \label{pr:trandeg}
Suppose that $p$ is a prime number, and let $d$, $e$, and
$\Upsilon\in\{0,1\}$ be cardinal numbers such that if $\Upsilon=0$
then $1\leq d$, if $p>2$ then $1\leq e$, and if $p=2$ and
$\Upsilon=1$ then $1\leq e$.  Then there exists a field $F$
containing $\Q(\xi_p)$ and a cyclic Galois extension $K$ of degree
$p$ over $F$ such that
\begin{equation*}
    e(K/F)=e, \quad d(K/F)=d, \quad \text{and }\quad
    \Upsilon(K/F) = \Upsilon.
\end{equation*}
Moreover
\begin{equation*}
    \trdeg (F/\Q)\leq 1+ \max\{e,d+1\}.
    \end{equation*}
In particular if $d, e\in\N\cup\{0\}$, there exists a Galois
cyclic extension $K/F$ of degree $p$ with prescribed invariants
$d$, $e$, and $\Upsilon$ of finite transcendence degree over its
prime field $\Q$.
\end{proposition}

\section{Acknowledgements} \label{S4}

We would like to thank the organizers of the MSRI programs on
Galois theory and the MSRI staff for giving us the opportunity to
meet and to begin our collaboration in the Fall of 1999. We are
also grateful to A.~Wadsworth for stimulating conversations. The
first author is very appreciative of the kind assistance of Ron
Hemphill, manager of Ivest Properties Limited (London, Canada),
for having provided excellent working conditions.

\end{document}